\documentclass[a4paper,10pt]{scrartcl}
\usepackage[utf8]{inputenc}
\usepackage{amsmath,amssymb,todonotes}
\usepackage{fullpage}
\newcommand{\calS}{\mathcal{S}}
\newcommand{\calE}{\mathcal{E}}
\newcommand{\del}{\partial}
\newcommand{\sprod}[1]{\langle#1\rangle}
\newcommand{\norm}[1]{\|#1\|}
\newcommand{\R}{\mathbb{R}}
\newcommand{\Rd}{\R^3}
\newcommand{\Rnn}{\R^{n\times n}}
\newcommand{\dV}{\,{\rm d}V}
\DeclareMathOperator{\tr}{tr}
\newcommand{\na}{\nabla}
\newcommand{\phii}{\varphi}
\newcommand{\Wbar}{\overline{W}}
\newcommand{\Rbar}{\overline{R}}

\newcommand{\Rtilde}{\widetilde{R}}
\newcommand{\Sz}{\mathbb{S}^2}
\newcommand{\dS}{\,\mathrm{dS}}
\newcommand{\dr}{{\,\rm dr}}
\newcommand{\so}{\mathfrak{so}}
\newcommand{\SO}{\mathrm{SO}}
\DeclareMathOperator{\skw}{skew}
\DeclareMathOperator{\sym}{sym}
\DeclareMathOperator{\Log}{Log}
\newcommand{\matr}[1]{\begin{pmatrix}#1\end{pmatrix}}
\newcommand{\col}{\colon}
\newcommand{\ny}{\nu}
\newcommand{\my}{\mu}
\newcommand{\tel}[1]{\frac{1}{#1}}
\usepackage[T1]{fontenc}
\newcommand{\changefont}[3]{\fontfamily{#1} \fontseries{#2} \fontshape{#3} \selectfont}

\newcommand{\orig}[1]{{\begin{quote} \changefont{ptm}{b}{it}
{\small \boldmath #1}\end{quote}}}
\title{On Grioli's minimum property and its relation to Cauchy's polar decomposition}
\author{Patrizio Neff\thanks{Lehrstuhl für Nichtlineare Analysis und Modellierung, Fakultät für Mathematik, Universität Duisburg-Essen, Thea-Leymann-Str. 9, 45127 Essen, Germany.}~\thanks{International Center M\&MOCS, Palazzo Caetani, Cisterna di Latina, Italy.}~\thanks{To whom correspondence should be addressed. e-mail: {\tt patrizio.neff@uni-due.de}}\quad and\quad  Johannes Lankeit\footnotemark[1]\quad and\quad Angela Madeo\thanks{Université de Lyon-INSA, 20 Av. Albert Einstein, 69100 Villeurbanne Cedex, France.}~\footnotemark[2]}
\date{}
\begin{document}
\maketitle
\begin{abstract}
  {{\bf Keywords:} polar decomposition, optimality of the polar factor, euclidean distance, geodesic distance, euclidean movement}
 \end{abstract}
 
 Every invertible matrix $F\in\Rnn$ can be uniquely decomposed into a product of a unitary matrix $R\in O(n)$ and a positive definite matrix $U$:
 \[
  F=R\,U.
 \]
 The roots of this ``polar decomposition theorem'' lie in Cauchy's work on elasticity \cite{Cauchy1841}. Finger gave it as an algebraic statement and ideas for a proof \cite[Eq (25)]{Finger1892}, the brothers E. and F. Cosserat proved it \cite[§6]{Cosserat1896}. Matrix notation and extension to the complex case are due to Autonne \cite{Autonne1902}, cf. \cite[Sect. 43]{Ericksen60}, \cite[Sect. 35-37]{Truesdell60}.
  (The result also holds for complex matrices and for non-square matrices (then upon loosing the uniqueness property of $R$), see e.g. \cite[ch. 8]{Higham2008}.)
  
 The unitary polar factor $R$ plays an important role in geometrically exact descriptions of solid materials. In this case  $R^TF=U$ is called the right stretch tensor of the deformation gradient $F=\nabla \phii$ and serves as a basic measure of the elastic deformation \cite{NeffBirsan13,Neff_Chelminski_ifb07,Neff_Forest_jel05,Neff_micromorphic_rse_05,Neff01c}. Indeed, it is known that the strain energy density for isotropic materials must depend only on the stretch $U$ in order to be frame-indifferent. Similar reasonings on objectivity lead to the result that the strain energy for isotropic second gradient materials must depend on the stretch $U$ and on its spatial gradient  (see \cite{IsolaGeneralizedHooke, EdgeIsolaSepp, EdgeContactIsolaSeppecher}). For additional applications and computational issues of the polar decomposition see e.g. \cite[Ch.~12]{Golubbookori} and \cite{nakahig12, ByXu08,nakatsukasa:2700,nahi11}.
 
The unitary polar factor can be characterized by its best-approximation property. For given $F$, it is the unique unitary matrix realizing
\begin{equation}
 \label{eq:min}
 \inf_{Q}\norm{F-Q}^2=\inf_Q\norm{Q^TF-I}^2=\norm{\sqrt{F^TF}-I}^2=\norm{U-I}^2
\end{equation}
over all unitary matrices $Q$, where $\norm{\cdot}$ is an arbitrary unitarily invariant norm \cite{Fan55}.

Optimality of the unitary polar factor is presently shown even for the expression $\norm{\Log Q^TF}$ in \cite{LankeitNeffNakatsukasa2013}, a distance measure arising from geometric considerations, connected with geodesic distances on matrix Lie groups (see \cite{NNF}, \cite{Neff_log_inequality13}, \cite{NeffEidelOsterbrinkMartin_Riemannianapproach} and \cite{Neff_Osterbrink_Martin_hencky13}). Here, $\Log$ is the (possibly multi-valued) matrix logarithm, i.e. a solution of $\exp(X)=Q^TF$.
In contrast to $\norm{F-Q}$ (cf. \cite{Neff_Biot07}), in this logarithmic expression symmetric and skew-symmetric part of the matrix norm can be weighted differently and the optimality of the polar factor still holds: 
\[\min_{Q\in \SO(n)}\left(\my\norm{\sym\Log(Q^TF)}^2+\my_c\norm{\skw\Log(Q^TF)}^2\right)= \mu\norm{\log \sqrt{F^TF}}^2\]
for $\my>0, \my_c\geq 0$, whereas the unitary polar factor fails to minimize the weighted expression
\[
 \my\norm{\sym(Q^TF-I)}^2+\my_c\norm{\skw(Q^TF-I)}^2, \quad 0<\my_c<\my
\]
in the Frobenius norm.

In this short note we would like to trace back the development on the optimality of the polar factor to its presumable roots, the work \cite{Grioli40} of G. Grioli, who shows the minimization property \eqref{eq:min} in the important special case of (some expression equivalent to) the Frobenius matrix norm and dimension $3$.

This work seems to have gone nearly unnoticed (but \cite{Guidugli79}, \cite{Guidugli80}, \cite{Pietraszkiewicz05} and \cite{Truesdell60}) and certainly the matrix-analysis community seems not to be aware of it. 
For example, \cite{Higham2008} refers to the work \cite{Fan55} of Fan and Hoffman for the optimality property (this is quite natural when being concerned with all unitarily invariant norms), who in turn seem to be nescient of Grioli's work.

We will juxtapose Grioli's original work \cite{Grioli40}, carefully translated from the original Italian paper by us, to a version with current notation. It will become clear that Grioli is showing even more: He considers Euclidean movements, not only linear transformations. Therefore, in his framework, it is not possible to consider weighted expressions.

While our paper does not contain new original results it may serve a pedagogical purpose: fundamental results are always older than it appears (see e.g. \cite{RussoForgotten}). 
\\[0.5cm]

Grioli starts by putting himself in the framework of finitely deforming bodies:

\orig{
Let $C_*$ be the reference configuration of an arbitrary continuous material system $\calS$; $C$ and $C'$ the current configurations of $\calS$ as a consequence of two different regular displacements $S$ and $S'$, $P_*$ the generic point of $C_*$; $P$ the corresponding of $P_*$.}
We consider a domain $C_*$, an arbitrary point $\vec p_*\in C_*$ and diffeomorphisms 
\begin{align*}
 S\col C_*\to C,\qquad\qquad S'\col C_*\to C'
\end{align*}
and denote $\vec p=S(\vec p_*)$, ${\vec p}\,'=S'(\vec p_*)$. We then restrict our investigation to a small ball $c_*=B_\rho(\vec p_*)$, where the affine approximation (by the first terms of the Taylor expansion)
\begin{equation}
\label{eq:approx}
 S(\vec p_*+h) \approx \vec p+\nabla S(\vec p_*).h
\end{equation}
is sufficiently good. 

\orig{Let then $c_*$ be a sphere with center $P_*$ and radius $\rho$ very small, which must be intended to be fixed independently of $P_*$. 

More precisely, we will consider $\rho$ to be so small that (correspondingly to any $P_*$) the displacements $S$ and $S'$ in $c_*$ can be identified with the corresponding homogeneous displacements tangent in $P_*$. If the displacements $S$ and $S'$ were homogeneous, no limitation would exist for $\rho$.

With reference to the arbitrary point $P_*$ it is common to define "local distance" of the two displacements $S$ and $S'$ the integral:
\[
d_{P_*}=\int_{c_*} |Q'Q|^2\, dC_*,
\]
where $Q$ and $Q'$ are the corresponding points in $C$ and $C'$ respectively to the arbitrary point $Q_*$ of $c_*$.}

The distance Grioli uses is 
\[
 d_{\vec p_*}(S,S')=\underset{x\in B_\rho(\vec p_*)}{\int} |S(x)-S'(x)|_{\R^3}^2\, dV=\underset{h\in B_\rho(0)}{\int}|S(\vec p_*+h)-S'(\vec p_*+h)|^2_{\Rd}\dV.
\]
To understand this distance measure better and demonstrate its connections to the Frobenius norm $\norm{Z}_F=\sqrt{\tr Z^TZ}$, for the moment we assume $Z$ and $Z'$ to be linear. Then 
\begin{align*}
 d_{\vec p_*}(Z,Z')&=\underset{x\in B_\rho(\vec p_*)}{\int}\sprod{Z(x)-Z'(x),Z(x)-Z'(x)}dx\\
 &=\underset{x\in B_\rho(\vec p_*)}{\int}\sprod{(Z-Z')^T(Z-Z')x,x}dx\\
 &=\frac{4\pi\rho^5}{15}\tr((Z-Z')^T(Z-Z'))
 =\frac{4\pi\rho^5}{15}\norm{Z-Z'}_F^2.
\end{align*}
Herein, the Frobenius-norm is obtained, since 
\begin{align}
\label{eq:Frob}
 \int_{h\in B_\rho(0)}\sprod{Zh,h}\dV&=\int_{h\in B_\rho(0)}\sum_{i,j=1}^3 z_{ij}h_ih_j \dV \notag
 \\&=\int_{h\in B_\rho(0)} z_{11}h_1^2+z_{22}h_2^2+z_{33}h_3^2 \dV = \tr Z\int_{h\in B_\rho(0)} h_1^2 \dV
 \\&=\tr Z\int_{h\in B_\rho(0)} \frac{h_1^2+h_2^2+h_3^2}{3} \dV =\tr Z \int_{r=0}^\rho \int_{\Sz} \frac{r^2}{3}  \dS\; r^2\, \dr
 =\tr Z \,\frac{4\pi\rho^5}{15},\notag
\end{align}
where the third equality holds since $\int_h h_i^2\dV = \int_h h_j^2\dV, i,j=1,2,3$ and where we used that $\int_{\Sz} 1 \dS=4\pi$.

The integration over the sphere/ball/... is a useful concept in order to average out (homogenize) direction dependent response. For applications in gradient elasticity, see e.g. \cite{Neff_Jeong_IJSS09, Neff_Jeong_ZAMP08}. 
It leads in a natural way to an isotropization, at the expense, hovever, to oversimplify the material response in special cases. This is e.g. the case in linear elasticity theory 
where it leads to a Poisson number $\ny=\tel4$:

Consider a linear elastic body, the strain energy of a small homogeneous sample in response to a displacement $u$ is to be obtained.
Locally, the energy should be quadratic in the distance of neighbouring particles.
Let $x$ and $x+h$ be two such particles. The elastic interaction in the direction $h$ is governed by a quadratic spring with spring constant $\my>0$. Hence, the directional energy is 
\begin{equation}
\label{eq:energy}
 \calE_h(x)= \frac \my 2 \sprod{u(x+h)-u(x),h}_{\R^3}^2.
\end{equation}
Since no direction is preferred in an isotropic body, the dependence on the direction can be averaged out and the total energy is obtained as integral over a sphere 
\[
 \calE(x)=
 \frac\my2\int_{h\in\Sz}\calE_{h}(x)\dS.
\]
Assuming a Taylor expansion $u(x+h)=u(x)+\nabla u(x).h+\ldots$ and approximating \eqref{eq:energy} by $\sprod{\nabla u(x).h,h}^2$, i.e.
\[
 \calE(x) \sim
 \frac\my2\int_{h\in \Sz} \sprod{\nabla u(x).h,h}^2 \dS
\]
and using (cf. \cite{Neff_Jeong_IJSS09})
 \[
  \int_{h\in\Sz}\sprod{Z.h,h}^2\dS=\frac{4\pi}{15}\left(2\norm{\sym Z}^2+[\tr Z]^2 \right),
 \]
one arrives at 
\[
 \calE(x)=
 \frac{4\pi}{15}\big(\my\norm{\sym \nabla u(x)}^2+\underbrace{\frac\my2}_{=:\frac\lambda2}[\tr \nabla u(x)]^2 \big),
\]
which corresponds to $\ny=\frac\lambda{2(\my+\lambda)}=\tel4$.

\orig{If one thinks the regular displacement $S$ to be arbitrarily assigned, one can ask himself: corresponding to arbitrary $P_*$, what is the rigid displacement which has the minimum local distance from $S$? What is, in other words, the rigid displacement which, locally, best approximates $S$?}
Grioli aims to find a rigid $S_{r(igid)}'$, such that $d_{\vec p_*}(S,S_r')$ is minimal. ``Rigid displacement'' means that $S_r'$ is of the form $S_r'(x)=t'+R'x$ for some $t'\in\R^3$ and $R'\in \SO(3)$, that is a constant rotation followed by a constant translation:
\[
 \inf_{t'\in \R^3, \Rbar\in \SO(3)} \int_{h\in B_\rho(0)}|S(\vec p_*+h)-[\Rbar.h+t']|_{\Rd}^2\dV.
\]
This problem becomes simpler if not rotations, but their infinitesimal version, represented by elements of $\so(3)$, i.e. skew-symmetric matrices, are considered, and Grioli gives a reference: 
\orig{For an infinitesimal displacement $S$ the answer has already been known for a long time\footnote{\changefont{ptm}{b}{n} Sobrero, Lezioni di Fisica Matematica, Roma, 1935-36}: decomposing the displacement (homogeneous, infinitesimal) in a rigid displacement plus a pure deformation, one gets as rigid displacement exactly the one which best approximates the effective displacement of the particle.}
For $W\in \so(3)$ and for an arbitrary matrix $\nabla u$ (in linear elasticity theory, usually the displacement gradient $\nabla u$ is of interest)
\begin{align*}
 \inf_{t'\in\Rd, \Wbar\in \so(3)} &\int_{h\in B_\rho(0)} |\vec p+\na u.h-[\Wbar.h+ t']|_{\Rd}^2\dV \\
 &= \inf_{t'\in\Rd, \Wbar\in \so(3)} \int_{h\in B_\rho(0)} |\vec p-t'|_{\Rd}^2 + |(\nabla u-\Wbar).h|_{\Rd}^2 + 2 \underbrace{\sprod{\vec p-t',(\nabla u-\Wbar).h}}_{=0,\textrm{ after integration (symmetry)}} \dV\\
 &\overset{\textrm{choose } t':=\vec p}{=}\inf_{\Wbar \in\so(3)} \int_{h\in B_\rho(0)}|(\nabla u-\Wbar).h|^2_{\Rd}\dV\\&= \frac{4\pi\rho^5}{15}\inf_{\Wbar \in\so(3)} \norm{\nabla u-\Wbar}^2_F\\ &= \frac{4\pi\rho^5}{15}\norm{\sym \na u}^2_F,\qquad \Wbar=\skw \na u.
\end{align*}
This development shows one way which allows to motivate the small strain tensor $\varepsilon=\sym\nabla u$ of linearized elasticity theory.

\orig{We want to show how an analogous theorem exists also if the deformation $S$ is not infinitesimal, but which is based on the requirement that the homogeneous displacement tangent to $C$ in $P_*$ is decomposed (which is of course possible) into the product of a pure deformation $D^*$ and a rigid displacement $S_r^*$.}
Using the left polar decomposition to express $\nabla S$ as a product of a rotation $R\in SO(3)$ and a pure deformation $D$ with eigenvalues $1+\Delta_1, 1+\Delta_2, 1+\Delta_3$, that is
\[
 \nabla S=R\,D,
\]
Grioli sets out to deduce a lower bound for $d_{\vec p_*}(S,S'_r)$ in terms of $\Delta_i$, that is, in terms of the positive definite polar factor of $\nabla S$.
\orig{More precisely, if we indicate by $\Delta_1, \Delta_2, \Delta_3$ the principal coefficients of the linear dilatation in $P_*$, we will show that (with reference to $P_*$) the local difference of $S$ to {\it any} rigid displacement $S_r'$ is always such that
\[
 d_{P_*}\geq\frac{4}{15}\pi\rho^5\,(\Delta_1^2+\Delta_2^2+\Delta_3^2),
\]
where equality holds if and only if $S_r'$ coincides with $S_r^*$.}
\[
 \inf_{t'\in \R^3, \Rbar\in \SO(3)} \int_{h\in B_\rho(0)}|S(\vec p_*+h)-[\Rbar.h+t']|^2\dV=\frac{4}{15}\pi\rho^5\,(\Delta_1^2+\Delta_2^2+\Delta_3^2)
\]

\orig{Also in the proof we will keep using the notations used by prof. SIGNORINI in different publications\footnote{\changefont{ptm}{b}{n} A. Signorini: Atti Accad. Naz. Lincei. Rendiconti 1930, vol. XII, p. 312:{\changefont{ptm}{b}{it} Sulle deformazioni finite dei sistemi continui}.} and in his current course on finite elastic transformations held at the National Institute of High Mathematics.
\begin{center}
 ***
\end{center}}

Grioli starts his proof by choosing an appropriate coordinate system: the coordinate axes are the eigenvectors of the positive definite part in the polar decomposition of $\nabla S=RD$.
\orig{With reference to a specific but arbitrary point $P_*$ we choose the Cartesian reference frame $T=P_* i_1, i_2, i_3$ with the condition that it is a principal system for deformation (in $P_*$) and, with respect to this system, we denote by $y_1, y_2, y_3$ the coordinates of the generic $Q_*$  and by $x_1, x_2, x_3$ those of the corresponding $Q$.

We also set
\[
 \alpha\equiv\left\|\begin{matrix}\frac{\del x_1}{\del y_1}&\frac{\del x_1}{\del y_2}&\frac{\del x_1}{\del y_3}\\\frac{\del x_2}{\del y_1}&\frac{\del x_2}{\del y_2}&\frac{\del x_2}{\del y_3}\\\frac{\del x_3}{\del y_1}&\frac{\del x_3}{\del y_2}&\frac{\del x_3}{\del y_3}\end{matrix}\right\|,
\]
where all the $\frac{\del x_r}{\del y_s} (r,s=1,2,3)$ are meant to be calculated in $P_*$.}
This mapping $\alpha$ corresponds to $\nabla S(\vec p_*)$ in our notation. As already announced, he then decomposes 
\[
 \nabla S(\vec p_*) =R\,D:
\]

\orig{By virtue of the regularity conditions for the displacement $S$, the $\alpha$ can be of course decomposed in the left product of a dilatation $\alpha_d$ with principal coefficients which are all positive and a rotation, $\alpha_r$: $\alpha_d$ characterizes (the dilatation)  $D$ while $S_r^*$ is the product of the rigid rotation characterized by $\alpha_r$ and the translation $P_*P$.

In other terms we have 
\[
 P_*Q=P_*P+PQ=P_*P+\alpha_r\alpha_d(P_*Q_*),
\]}
that is
\[
 S(x)-\vec p_*=S(\vec p_*)-\vec p_* +\nabla S(\vec p_*)(x-\vec p_*)= \vec p-\vec p_* + RD(x-\vec p_*)
\]
(Of course, the first equality sign uses hypothesis \eqref{eq:approx}.)

The substraction of $\vec p_*$ is due to his notation for vectors as vectors from one point to another; we would maybe rather write the equivalent 
\[
 S(x)=S(\vec p_*)+\nabla S(\vec p_*)(x-\vec p_*)= \vec p+ R\,D(x-\vec p_*).
\]

Also for the wanted rigid movement such an expression can be given: a rotation $R'$ centered in $\vec p_*$ followed by a translation by some $t'$:
\[
 S'_r(x)-\vec p_*=t'+R'(x-\vec p_*) ...
\]

\orig{while if we consider also $S_r'$ as the product of a rigid rotation characterized by $\alpha_r'$ and of a translation $t'$ we can set
\[
P_*Q'=t'+\alpha_r'(P_*Q_*) 
\]
and hence
\begin{equation}
 \label{eq:one}
 Q'Q=P_*P-t'+\alpha_r\alpha_d(P_*Q_*)-\alpha_r'(P_*Q_*).
\end{equation}}
... which leads to this representation of the difference between $S$ and $S'_r$. 
\[
 S(x)-S'_r(x)=\vec p-\vec p_*-t' +RD(x-\vec p_*)-R'(x-\vec p_*).
\]
$D$ is positive definite and, according to the choice of the coordinate system, diagonal:
\begin{equation}
\label{eq:defD}
 D=\matr{1+\Delta_1&&\\&1+\Delta_2&\\&&1+\Delta_3}.
\end{equation}

\orig{Moreover, for the way in which we chose $T$,  we have 
\begin{equation}
 \label{eq:two} 
 \alpha_d\equiv\left\|\begin{matrix}
                       1+\Delta_1&0&0\\0&1+\Delta_2&0\\0&0&1+\Delta_3
                      \end{matrix}\right\|
\end{equation}
with
\begin{equation}
 \label{eq:twoprime}
 1+\Delta_r>0, \qquad (r=1,2,3)\tag{2'} 
\end{equation}
and also 
\begin{align}
 \label{eq:three}
 \int_{c_*}y_r\,dC_*=0,&\quad \int_{c_*}y_r\,y_s\,dC_*=0,\quad\int_{c_*}y_r^2\,dC_*=\frac{4}{15}\pi\rho^5,\\ &\qquad(r,s=1,2,3; r\neq s).&\notag
\end{align}

\begin{center}
 ***
\end{center}}
\orig{Since the local distance only depends on the length of the vectors $Q'Q$, $d_{P_*}$  is not changed if we replace $Q'Q$ with $\alpha_r^{-1}(Q'Q)$, since the invertible linear mapping $\alpha_r^{-1}$ is also a rotation and hence does not affect the distances.}
Grioli then notes that $d_{\vec p_*}(S,S'_r)=d_{\vec p_*}(R^TS,R^TS_r')$, since $|y|_{\Rd}=|R^Ty|_{\Rd}$ and therefore obtains
\[
 d_{\vec p_*}(S,S_r')=d_{\vec p_*}(R^TS,R^TS_r')=\underset{B_\rho(\vec p_*)}{\int}|\underbrace{R^T(\vec p-\vec p_*-t')}_{=:t''} +D(x-\vec p_*)-\underbrace{R^TR'}_{=:\Rtilde}(x-\vec p_*)|^2_{\Rd}\dV
\]
where he defines (what we will call) $t''$ and $\Rtilde$.

\orig{We set
\begin{equation}
 \label{eq:four}
 t''=\alpha_r^{-1}(P_*P-t'),\qquad \alpha_r''=\alpha_r^{-1}\alpha_r',
\end{equation}
so that also the invertible linear mapping $\alpha_r''$ turns out to be a rotation.

Using \eqref{eq:one} we get
\begin{equation}
 \label{eq:five}
 d_{P_*}=\int_{c_*}dC_*|t+\alpha_d(P_*Q_*)-\alpha_r''(P_*Q_*)|^2.
\end{equation}}
He then computes this integral $d_{\vec p_*}(S,S'_r)$.

\orig{Using \eqref{eq:two} and \eqref{eq:three} and the trivial equality
\[
 |\alpha_r''(P_*Q_*)|^2=|P_*Q_*|^2,
\]
equation \eqref{eq:five} is easily simplified in
\[
 d_{P_*}=\frac{4}{3}\pi\rho^3t''^2+\left|{\sum_{s=1}^3}(1+\Delta_s)^2+3\right|\;\frac{4}{15}\pi\rho^5-2\int_{c_*}dC_*\;\alpha_d(P_*Q_*)\cdot\alpha_r''(P_*Q_*).
\]}
Here he has used 
\begin{align*}
 &|t'' +D(x-\vec p_*)-\Rtilde(x-\vec p_*)|^2_{\Rd}\\&\quad=|t''|_{\Rd}^2+|D(x-\vec p_*)|_{\Rd}^2+|\Rtilde (x-\vec p_*)|_{\Rd}^2\\&\qquad +2\sprod{t'',(D-\Rtilde)(x-\vec p_*)}\;-\;2\sprod{D(x-\vec p_*),\Rtilde(x-\vec p_*)},
\end{align*}
where the first scalar product term vanishes upon integration (because of symmetry) and as rotation and therefore isometry $\Rtilde$ in $|\Rtilde (x-\vec p_*)|^2$ can be neglected.

What remains, is 
\begin{equation*}
 \underset{B_\rho(\vec p_*)}{\int} |t''|_{\Rd}^2 \dV+ \underset{B_\rho(\vec p_*)}{\int} |D(x-\vec p_*)|_{\Rd}^2 \dV+ \underset{B_\rho(\vec p_*)}{\int} |(x-\vec p_*)|_{\Rd}^2\dV - 2\underset{B_\rho(\vec p_*)}{\int} \sprod{D(x-\vec p_*),R(x-\vec p_*)}\dV.
\end{equation*}
Herein, 
\[
 \underset{B_\rho(\vec p_*)}{\int} |t''|_{\Rd}^2 \dV = \frac{4}{3}\pi\rho^3 |t''|_{\Rd}^2,
\]
and (confer \eqref{eq:Frob} and \eqref{eq:defD}) 
\[
 \underset{B_\rho(\vec p_*)}{\int} |D(x-\vec p_*)|^2_{\Rd} \dV=\norm{D}_F^2\frac{4\pi}{15}\rho^5 = \sum_{s=1}^3 (1+\Delta_s)^2 \frac{4\pi}{15}\rho^5. 
\]
Analogously, $\int_{B_\rho} |(x-\vec p_*)|^2=\norm{I}_F^2\frac{4\pi}{15}\rho^5=3\frac{4\pi}{15}\rho^5$. Next we compute the value of the integral over the scalar product term: 
\orig{On the other hand, if we call $c_{rs}$ the coefficients of the invertible linear mapping $\alpha_r''$ with respect to $T$, we have [cfr. again \eqref{eq:two} and \eqref{eq:three}]
\begin{align*}
 \int_{c_*}dC_*\alpha_d(P_*Q_*)\times \alpha_r''(P_*Q_*)=&\\
 =\int_{c_*}dC_*\left|{\sum_{s=1}^3} (1+\Delta_s)y_si_s\times{\sum_{r,s=1}^3}c_{rs}y_si_r\right|=&\frac{4}{15}\pi\rho^5{\sum_{s=1}^3}(1+\Delta_s)c_{ss},
\end{align*}}

We let $\Rtilde=\matr{r_{11}&r_{12}&r_{13}\\r_{21}&r_{22}&r_{23}\\r_{31}&r_{32}&r_{33}}$, where $|r_{ij}|\leq 1$, because $\Rtilde$ is orthogonal, ($r_{ij}=c_{ij}$ in Grioli's terminology) and observe
\[
 - 2\underset{B_\rho(\vec p_*)}{\int} \sprod{D(x-\vec p_*),\Rtilde(x-\vec p_*)}\dV = - 2\underset{h\in B_\rho(0)}{\int} \sprod{\matr{(1+\Delta_1)h_1\\(1+\Delta_2)h_2\\(1+\Delta_3)h_3},\matr{r_{11}h_1+r_{12}h_2+r_{13}h_3\\r_{21}h_1+r_{22}h_2+r_{23}h_3\\r_{31}h_1+r_{32}h_2+r_{33}h_3}}\dV.
\]
As Grioli has computed earlier (cf. \eqref{eq:four}), 
the mixed terms $h_1h_2$ etc. yield $0$. We are left with 
\[
 - 2 \underset{h\in B_\rho(0)}{\int} (1+\Delta_1)r_{11}h_1^2+(1+\Delta_2)r_{22}h_2^2+(1+\Delta_3)r_{33}h_3^2 \dV= \frac{4\pi}{15}\rho^5 \sum_{s=1}^3(-2(1+\Delta_s)r_{ss}).
\]

Combining all these, up to now it is shown that 
\[
 d_{\vec p_*}(S,S_r') = \frac{4}{3}\pi\rho^3 |t''|^2 +\frac{4\pi}{15}\rho^5 \left( \sum_{s=1}^3(1+\Delta_s)^2 +3 -2\sum_{s=1}^3(1+\Delta_s)r_{ss}\right).
\]
\orig{and hence we have the expression
\begin{equation}
 \label{eq:six}
 d_{P_*}=\frac{4}{3}\pi\rho^3t''^2+\frac{4}{15}\pi\rho^5|{\sum_{s=1}^3}(1+\Delta_s)^2+3-2(1+\Delta_s)c_{ss}|. 
\end{equation}}

This expression obviously becomes minimal if $t''$ vanishes ($t'=\vec p-\vec p_*$) and $r_{ss}$ is maximal - that is, if $r_{ss}=1$.
The conditions $r_{11}=r_{22}=r_{33}=1$ is not only, as Grioli remarks, contained in but equivalent to the condition that $\Rtilde$ is the identity. ($\Rtilde$ is orthogonal, so each column has to be a unit vector.)

The minimal value attained then is 
\begin{align*}
 d_{\vec p_*}(S,S_r') &= \frac{4}{3}\pi\rho^3 |0|^2 +\frac{4\pi}{15}\rho^5 \left( \sum_{s=1}^3\left((1+\Delta_s)^2 -2 (1+\Delta_s)1\right) + 3\right)\\
 &= \frac{4\pi}{15}\rho^5 (\Delta_1^2+\Delta_2^2+\Delta_3^2)=\frac{4\pi}{15}\rho^5\norm{D-I}_F^2.
\end{align*}

\orig{Let us now consider \eqref{eq:twoprime} and the fact that, since  $c_{rs}$ are all direction cosines it must be
\[
 c_{ss}\leq 1\qquad (s=1,2,3).
\]
This is sufficient to deduce from \eqref{eq:six}
\[
 d_{P_*}\geq\frac{4}{15}\pi\rho^5\,\left|{\sum_{s=1}^3}(1+\Delta_s)^2+3-2(1+\Delta_s)\right|,
\]
and hence
\begin{equation}
 \label{eq:seven}
  d_{P_*}\geq\frac{4}{15}\pi\rho^5|\Delta_1^2+\Delta_2^2+\Delta_3^2|.
\end{equation}
The equality sign holds if and only if one simultaneously has
\[
 t''=0,\qquad c_{ss}=1, \qquad (s=1,2,3).
\]
It is evident that the three conditions $c_{ss}=1 (s=1,2,3)$ are contained in the condition according to which the rotation $\alpha_r''$ reduces to the identity.

Hence, simply recalling \eqref{eq:four}, we can conclude that in \eqref{eq:seven} the equality sign holds if and only if one has simultaneously 
\[
 t=P_*P,\qquad \alpha_r'=\alpha_r,
\]
i.e. if and only if $S_r'$ coincides with $S_r^*$, qed.}
Thus, Grioli has shown that 
\begin{align*}
  \inf_{t'\in \R^3, \Rbar\in \SO(3)} \int_{h\in B_\rho(0)}|S(\vec p_*)&+\nabla S(\vec p_*).h -[\Rbar.h+t']|_{\Rd}^2\dV\\&=\frac{4\pi\rho^5}{15}\norm{D-I}^2_F
  =\frac{4\pi\rho^5}{15}\norm{\sqrt{\nabla S^T\nabla S}-I}^2_F,
\end{align*}
where $D$ is unitarily similar to the Hermitian polar factor $\sqrt{\nabla S^T\nabla S}$ of $\nabla S(\vec p_*)$.

{\footnotesize
    \bibliographystyle{plain} 

    \bibliography{literatur1}

\begin{thebibliography}{10}

\bibitem{Autonne1902}
L.~Autonne.
\newblock Sur les groupes lin{\'e}aires, r{\'e}els et orthogonaux.
\newblock {\em Bull. Soc. Math. France}, 30:121--134, 1902.

\bibitem{NeffBirsan13}
M.~B\^{i}rsan and P.~Neff.
\newblock Existence of minimizers in the geometrically non-linear 6-parameter
  resultant shell theory with drilling rotations.
\newblock {\em Mathematics and Mechanics of Solids}, 2013.

\bibitem{Neff_log_inequality13}
M.~B\^{i}rsan, P.~Neff, and J.~Lankeit.
\newblock Sum of squared logarithms: {A}n inequality relating positive definite
  matrices and their matrix logarithm.
\newblock {\em Journal of Inequalities and Applications}, 2013(1):168, 2013.

\bibitem{Pietraszkiewicz05}
C.~Bouby, D.~Fortun{\'e}, W.~Pietraszkiewicz, and C.~Vall{\'e}e.
\newblock Direct determination of the rotation in the polar decomposition of
  the deformation gradient by maximizing a {Ra}yleigh quotient.
\newblock {\em Z. Angew. Math. Mech.}, 85:155--162, 2005.

\bibitem{ByXu08}
R.~Byers and H.~Xu.
\newblock {A new scaling for Newton's iteration for the polar decomposition and
  its backward stability}.
\newblock {\em SIAM J. Matrix Anal. Appl.}, {30}:{822--843}, {2008}.

\bibitem{Cauchy1841}
A.~L. Cauchy.
\newblock {\em M\'emoire sur les dilatations, les condensations et les
  rotations produites par un changement de forme dans un syst\'eme de points
  mat\'eriel. Oeuvres compl\'etes d'Augustin Cauchy, II-XII, pp. 343--377}.
\newblock Gauthier-Villars, Paris, Paris, 1841.

\bibitem{Cosserat1896}
E.~Cosserat and F.~Cosserat.
\newblock Sur la th\'eorie de l'\'elasticit\'e.
\newblock {\em Ann. Toulouse}, 10:1--116, 1896.

\bibitem{IsolaGeneralizedHooke}
F.~dell'Isola, G.~Sciarra, and S.~Vidoli.
\newblock Generalized {H}ooke's law for isotropic second gradient materials.
\newblock {\em Proc. R. Soc. A}, 465(2107):2177--2196, 2009.

\bibitem{EdgeContactIsolaSeppecher}
F.~dell'Isola and P.~Seppecher.
\newblock The relationship between edge contact forces and intersticial working
  allowed by the principle of virtual power.
\newblock {\em Comptes Rendus M{\'e}canique}, t. 321(s{\'e}rie IIb):303--308,
  1995.

\bibitem{EdgeIsolaSepp}
F.~dell'Isola and P.~Seppecher.
\newblock Edge contact forces and quasi-balanced power.
\newblock {\em Meccanica}, 32(1):33--52, 1997.

\bibitem{Ericksen60}
J.L. Ericksen.
\newblock Tensor fields.
\newblock In S.~Fl\"ugge, editor, {\em Handbuch der {P}hysik}, volume III/1.
  Springer, Heidelberg, 1960.

\bibitem{Fan55}
K.~Fan and A.~J. Hoffmann.
\newblock Some metric inequalities in the space of matrices.
\newblock {\em Proc. Amer. Math. Soc.}, 6:111--116, 1955.

\bibitem{Finger1892}
J.~Finger.
\newblock {\"U}ber die gegenseitigen {B}eziehungen gewisser in der {M}echanik
  mit {V}ortheil anwendbaren {Fl\"a}chen zweiter {O}rdnung nebst {A}nwendungen
  auf {P}robleme der {A}statik.
\newblock {\em Wien Sitzgsber. 2a}, 101:1105--1142, 1892.

\bibitem{Golubbookori}
G.~H. Golub and C.~V.~Van Loan.
\newblock {\em Matrix Computations}.
\newblock {The Johns Hopkins University Press}, 1996.

\bibitem{Grioli40}
G.~Grioli.
\newblock Una propriet\`a di minimo nella cinematica delle deformazioni finite.
\newblock {\em Boll. Un. Math. Ital.}, 2:252--255, 1940.

\bibitem{Higham2008}
N.~J. Higham.
\newblock {\em Functions of {M}atrices: Theory and {C}omputation}.
\newblock SIAM, Philadelphia, PA, USA, 2008.

\bibitem{LankeitNeffNakatsukasa2013}
J.~{Lankeit}, P.~{Neff}, and Y.~{Nakatsukasa}.
\newblock {The minimization of matrix logarithms - {O}n a fundamental property
  of the unitary polar factor}.
\newblock {\em ArXiv e-prints}, August 2013.

\bibitem{Guidugli80}
L.C. Martins and P.~Podio-Guidugli.
\newblock An elementary proof of the polar decomposition theorem.
\newblock {\em The American Mathematical Monthly}, 87:288--290, 1980.

\bibitem{nakatsukasa:2700}
Y.~Nakatsukasa, Z.~Bai, and F.~Gygi.
\newblock {Optimizing Halley's iteration for computing the matrix polar
  decomposition}.
\newblock {\em SIAM J. Matrix Anal. Appl.}, 31(5):2700--2720, 2010.

\bibitem{nahi11}
Y.~Nakatsukasa and N.~J. Higham.
\newblock Backward stability of iterations for computing the polar
  decomposition.
\newblock {\em SIAM J. Matrix Anal. Appl.}, 33(2):460--479, 2012.

\bibitem{nakahig12}
Yuji Nakatsukasa and Nicholas~J. Higham.
\newblock Stable and efficient spectral divide and conquer algorithms for the
  symmetric eigenvalue decomposition and the svd.
\newblock {\em SIAM J. Sci. Comp}, 35(3):A1325--A1349, 2013.

\bibitem{Neff01c}
P.~Neff.
\newblock Local existence and uniqueness for quasistatic finite plasticity with
  grain boundary relaxation.
\newblock {\em Quart. Appl. Math.}, 63:88--116, 2005.

\bibitem{Neff_micromorphic_rse_05}
P.~Neff.
\newblock Existence of minimizers for a finite-strain micromorphic elastic
  solid.
\newblock {\em Proc. Roy. Soc. Edinb. A}, 136:997--1012, 2006.

\bibitem{Neff_Chelminski_ifb07}
P.~Neff and K.~Che{\l}mi\'nski.
\newblock A geometrically exact {C}osserat shell-model for defective elastic
  crystals. {J}ustification via {$\Gamma$}-convergence.
\newblock {\em Interfaces and Free Boundaries}, 9:455--492, 2007.

\bibitem{NeffEidelOsterbrinkMartin_Riemannianapproach}
P.~{Neff}, B.~{Eidel}, F.~{Osterbrink}, and R.~{Martin}.
\newblock {A Riemannian approach to strain measures in nonlinear elasticity}.
\newblock {\em ArXiv e-prints}, May 2013.
\newblock to appear in Comptes Rendus Mecanique.

\bibitem{Neff_Osterbrink_Martin_hencky13}
P.~Neff, B.~Eidel, F.~Osterbrink, and R.~Martin.
\newblock {T}he isotropic {H}encky strain energy measures the geodesic distance
  of the deformation gradient {$F \in\mathrm{GL^+}(n)$} to $\mathrm{SO}(n)$ in
  the unique left invariant {R}iemannian metric on $\mathrm{GL}(n)$ which is
  also right $\mathrm{O}(n)$-invariant.
\newblock {\em submitted}, 2013.

\bibitem{Neff_Biot07}
P.~Neff, A.~Fischle, and I.~M\"unch.
\newblock Symmetric {C}auchy-stresses do not imply symmetric {B}iot-strains in
  weak formulations of isotropic hyperelasticity with rotational degrees of
  freedom.
\newblock {\em Acta Mechanica}, 197:19--30, 2008.

\bibitem{Neff_Forest_jel05}
P.~Neff and S.~Forest.
\newblock A geometrically exact micromorphic model for elastic metallic foams
  accounting for affine microstructure. {M}odelling, existence of minimizers,
  identification of moduli and computational results.
\newblock {\em J. Elasticity}, 87:239--276, 2007.

\bibitem{Neff_Jeong_ZAMP08}
P.~Neff, J.~Jeong, I.~M\"unch, and H.~Ramezani.
\newblock Mean field modeling of isotropic random {C}auchy elasticity versus
  microstretch elasticity.
\newblock {\em Z. Angew. Math. Phys.}, 3(60):479--497, 2009.

\bibitem{Neff_Jeong_IJSS09}
P.~Neff, J.~Jeong, and H.~Ramezani.
\newblock Subgrid interaction and micro-randomness - novel invariance
  requirements in infinitesimal gradient elasticity.
\newblock {\em Int. J. Solids Struct.}, 46(25-26):4261--4276, 2009.

\bibitem{NNF}
P.~{Neff}, Y.~{Nakatsukasa}, and A.~{Fischle}.
\newblock The unitary polar factor ${Q}={U}_p$ minimizes $\|
  \rm{{L}og}$(${Q}^*${Z})$\|^2$ and $\|\rm{sym_*\, {L}og}$(${Q}^*${Z})$\|^2$ in
  the spectral norm in any dimension and the {F}robenius matrix norm in three
  dimensions.
\newblock {\em to appear in SIMAX}, 2013.
\newblock http://arxiv.org/abs/1302.3235.

\bibitem{Guidugli79}
P.~Podio-Guidugli and L.C. Martins.
\newblock A variational approach to the polar decomposition theorem.
\newblock {\em Rendiconti/Accademia Nazionale dei Lincei, Roma, Classe di
  Scienze Fisiche, Matematiche e Naturali}, LXVI (1979):487--493, 1979.

\bibitem{RussoForgotten}
L.~Russo.
\newblock {\em The forgotten revolution. How science was born in 300 BC and why
  it had to be reborn}.
\newblock Springer, 2004.

\bibitem{Truesdell60}
C.~Truesdell and R.~Toupin.
\newblock The classical field theories.
\newblock In S.~Fl\"ugge, editor, {\em Handbuch der {P}hysik}, volume III/1.
  Springer, Heidelberg, 1960.

\end{thebibliography}
}    

\section{Appendix}
 Let $Y=(Z-Z')^{T}(Z-Z')$, then
\begin{gather*}
d_{p*}\left(Z,Z'\right)=\int_{B_{\rho\left(p*\right)}}Y_{ij}x_{j}x_{i}\\
\\
Y_{11}\int_{B_{\rho\left(p*\right)}}x_{1}x_{1}+Y_{22}\int_{B_{\rho\left(p*\right)}}x_{2}x_{2}+Y_{33}\int_{B_{\rho\left(p*\right)}}x_{3}x_{3}\\
\\
+Y_{12}\int_{B_{\rho\left(p*\right)}}x_{1}x_{2}+Y_{13}\int_{B_{\rho\left(p*\right)}}x_{1}x_{3}+Y_{21}\int_{B_{\rho\left(p*\right)}}x_{2}x_{1}\\
\\
+Y_{23}\int_{B_{\rho\left(p*\right)}}x_{2}x_{3}+Y_{31}\int_{B_{\rho\left(p*\right)}}x_{3}x_{1}+Y_{32}\int_{B_{\rho\left(p*\right)}}x_{3}x_{2}
\end{gather*}

We know that to pass from cartesian to spherical coordinates one must
set
\begin{gather*}
x_{1}=\text{r sin\ensuremath{\theta}}\mathrm{cos}\phi\\
x_{2}=\text{r\ sin\ensuremath{\theta}}\mathrm{sin}\phi\\
x_{3}=\text{r\ cos\ensuremath{\theta}}
\end{gather*}
and that the Jacobian of the transformation is $J=r^{2}\:\text{sin\ensuremath{\theta}}$,
so that
\begin{gather*}
\int_{B_{\rho\left(p*\right)}}x_{1}x_{1}=\int_{0}^{\rho}r^{4}dr\int_{0}^{2\pi}\mathrm{cos}^{2}\phi\: d\phi\int_{0}^{\pi}\mathrm{sin}{}^{3}\text{\ensuremath{\theta}}\: d\theta=\frac{\rho^{5}}{5}\,\pi\,\frac{4}{3}=\frac{4\pi\rho^{5}}{15}\\
\\
\int_{B_{\rho\left(p*\right)}}x_{2}x_{2}=\int_{0}^{\rho}r^{4}dr\int_{0}^{2\pi}\mathrm{sin}^{2}\phi\: d\phi\int_{0}^{\pi}\mathrm{sin}{}^{3}\text{\ensuremath{\theta}}\: d\theta=\frac{\rho^{5}}{5}\,\pi\,\frac{4}{3}=\frac{4\pi\rho^{5}}{15}\\
\\
\int_{B_{\rho\left(p*\right)}}x_{3}x_{3}=\int_{0}^{\rho}r^{4}dr\int_{0}^{2\pi}\: d\phi\int_{0}^{\pi}\mathrm{cos}^{2}\theta\:\mathrm{sin}\text{\ensuremath{\theta}}\: d\theta=\frac{\rho^{5}}{5}\,2\pi\,\frac{2}{3}=\frac{4\pi\rho^{5}}{15}\\
\\
\int_{B_{\rho\left(p*\right)}}x_{1}x_{2}=\int_{B_{\rho\left(p*\right)}}x_{2}x_{1}=\int_{0}^{\rho}r^{4}dr\int_{0}^{2\pi}\mathrm{sin}\,\phi\mathrm{cos}\phi\: d\phi\int_{0}^{\pi}\mathrm{sin}{}^{3}\text{\ensuremath{\theta}}\: d\theta=\frac{\rho^{5}}{5}\,0\,\frac{4}{3}=0\\
\\
\int_{B_{\rho\left(p*\right)}}x_{1}x_{3}=\int_{B_{\rho\left(p*\right)}}x_{3}x_{1}=\int_{0}^{\rho}r^{4}dr\int_{0}^{2\pi}\mathrm{cos}\phi\: d\phi\int_{0}^{0}\mathrm{sin}{}^{2}\text{\ensuremath{\theta}}\mathrm{cos}\theta\: d\theta=\frac{\rho^{5}}{5}\,0\,0=0.\\
\\
\int_{B_{\rho\left(p*\right)}}x_{2}x_{3}=\int_{0}^{\rho}r^{4}dr\int_{0}^{2\pi}\mathrm{sin}\phi\: d\phi\int_{0}^{\pi}\mathrm{sin}^{2}\text{\ensuremath{\theta}}\mathrm{cos}\theta\: d\theta=\frac{\rho^{5}}{5}\,0\,0=0
\end{gather*}
Hence, we finally have
\begin{gather*}
d_{p*}\left(Z,Z'\right)=\left(Y_{11}+Y_{22}+Y_{33}\right)\frac{4\pi\rho^{5}}{15}
\end{gather*}
which is Eq. at the end of pag. 2.\\
A short version of this calculation can be found in \eqref{eq:Frob}.
\end{document}